\documentclass[12pt]{article}

\usepackage{amssymb}
\usepackage{amsmath,amsthm}
\usepackage[dvips]{graphicx}
\usepackage{wrapfig}

\newtheorem{thm}{Theorem}[section]
\newtheorem{cor}[thm]{Corollary}

\theoremstyle{definition}

\theoremstyle{remark}
\newtheorem{rem}[thm]{Remark.}

\title{Quantization of the crossing number of a knot diagram}
\author{Akio Kawauchi \thanks{Osaka City University Advanced Mathematical Institute, 3-3-138, Sugimoto, Sumiyoshi-ku Osaka 558-8585, Japan. kawauchi@sci.osaka-cu.ac.jp}
\and Ayaka Shimizu\thanks{Osaka City University Advanced Mathematical Institute, 3-3-138, Sugimoto, Sumiyoshi-ku Osaka 558-8585, Japan. shimizu1984@gmail.com}}
\date{\today}

\begin{document}

\maketitle

\begin{abstract}
We introduce the warping crossing polynomial of an oriented knot diagram by using the warping degrees of crossing points of the diagram. 
Given a closed transversely intersected plane curve, we consider oriented knot diagrams obtained from the plane curve as states  
to take the sum of the warping crossing polynomials for all the states for the plane curve. 
As an application, we show that every closed transversely intersected plane curve with even crossing points has two independent canonical orientations and every based closed transversely intersected plane curve with odd crossing points has two independent canonical orientations. 
\end{abstract}

\section{Introduction}

\begin{figure}[htbp]
 \begin{center}
  \includegraphics[width=30mm]{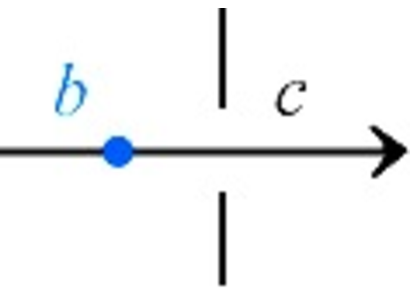}
 \end{center}
 \caption{}
 \label{c-wei}
\end{figure}

\begin{figure}[htbp]
 \begin{center}
  \includegraphics[width=30mm]{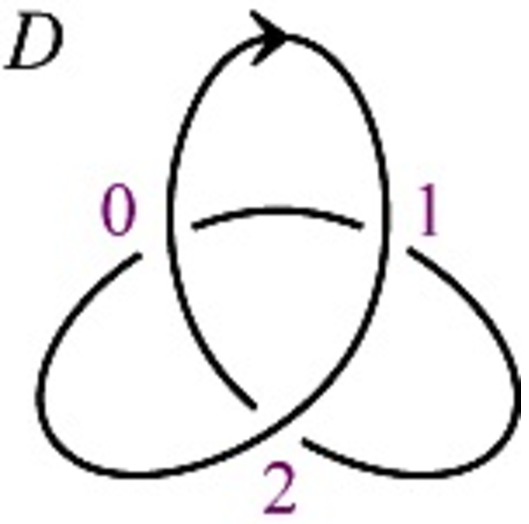}
 \end{center}
 \caption{}
 \label{x-ex}
\end{figure}

\begin{figure}[htbp]
 \begin{center}
  \includegraphics[width=120mm]{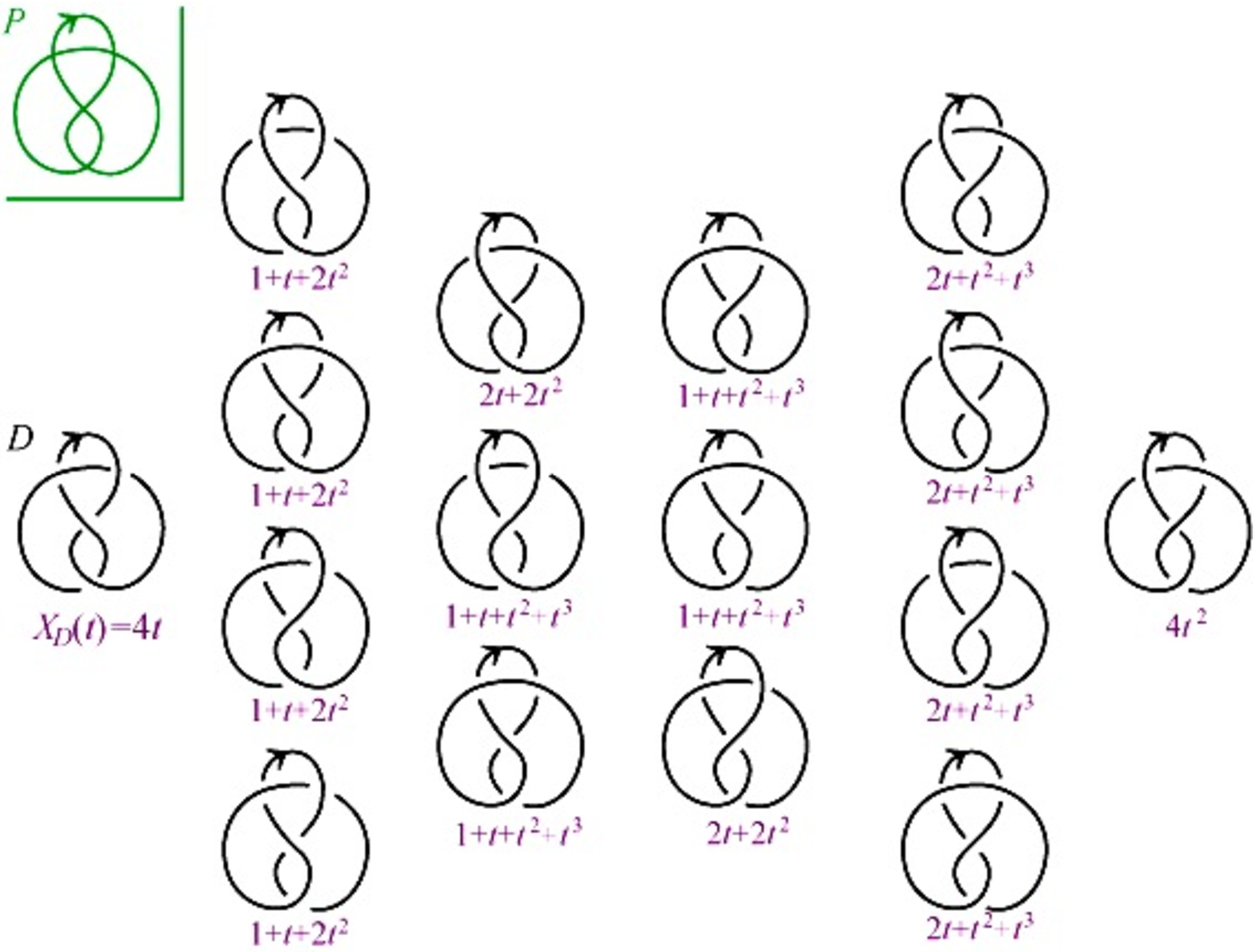}
 \end{center}
 \caption{}
 \label{s-sum1}
\end{figure}

\begin{figure}[htbp]
 \begin{center}
  \includegraphics[width=120mm]{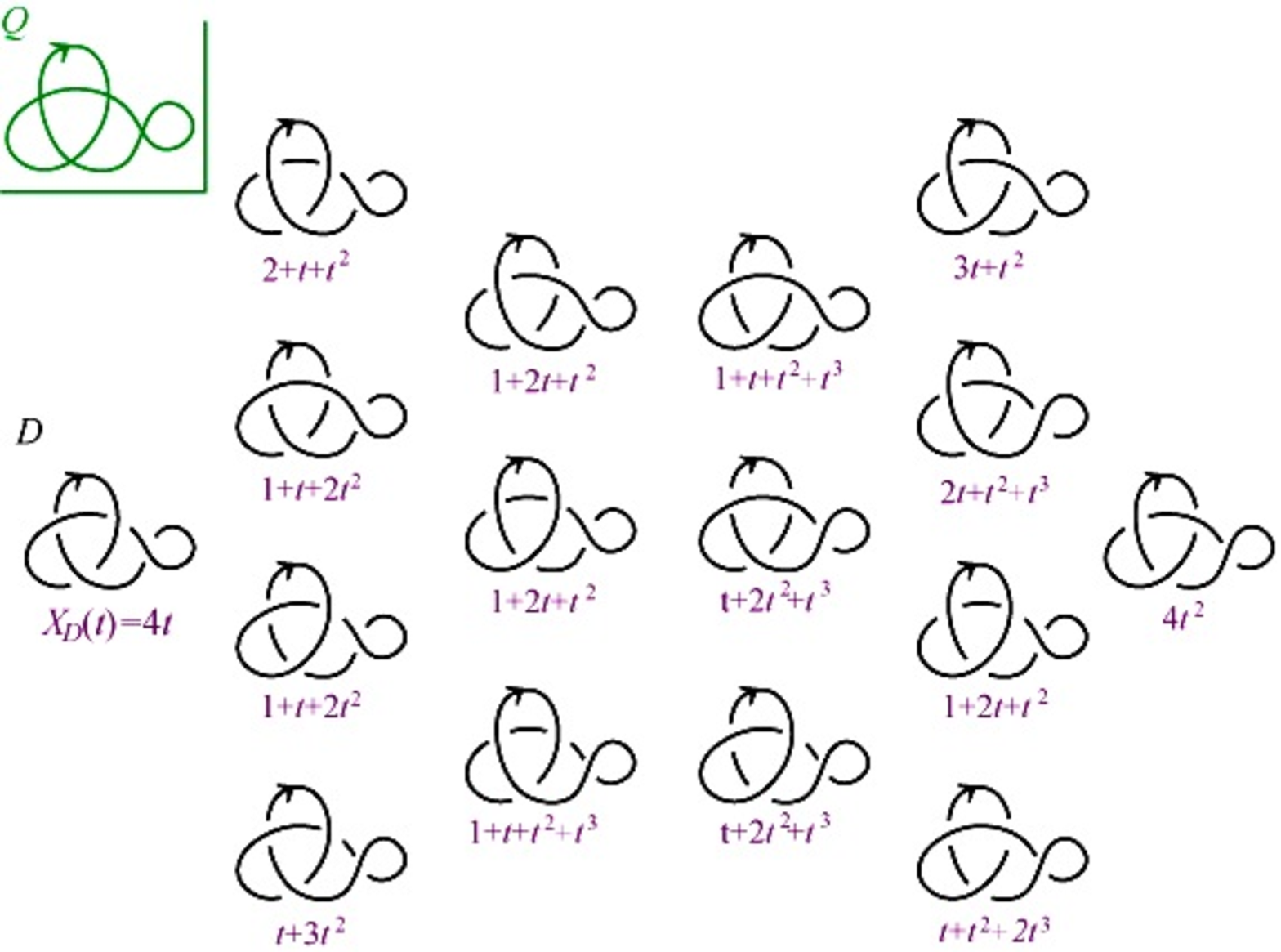}
 \end{center}
 \caption{}
 \label{s-sum2}
\end{figure}

Throughout this paper except Section \ref{s-ori}, knot diagrams are oriented and on $S^2$. 
A \textit{based diagram} $D_b$ is a diagram $D$ with a base point $b$. 
A crossing point of $D$ is a \textit{warping crossing point} of $D_b$ if we come to the crossing point as an under-crossing first when we go along $D$ with the orientation by starting from $b$. 
The \textit{warping degree} $d(D_b)$ of $D_b$ is the number of warping crossing points of $D_b$ \cite{kawauchi-l}. 
The warping degree is also defined for link diagrams and spatial graphs \cite{kawauchi-c}. 
We note that the similar notions are studied by Fujimura \cite{fujimura}, Fung \cite{fung}, Lickorish and Millett \cite{LM}, Okuda \cite{okuda} and Ozawa \cite{ozawa} considering the ascending number with an orientation. 
We define a weight of each crossing point $c$ of a knot diagram $D$ as follows: 
Take a base point $b$ which is just before the over-crossing of $c$ (Figure \ref{c-wei}). 
The \textit{crossing weight} $X_c(t)$ of $c$ is defined to be $t^{d(c)}$, where $d(c)=d(D_b)$. 
Now we define the \textit{warping crossing polynomial} $X_D(t)$ of a knot diagram $D$ to be 
the sum of crossing weights for all crossing points of $D$, i.e., $X_D(t)=\sum _c X_c(t)$. 
For example, the diagram $D$ in Figure \ref{x-ex} has $X_D(t)=1+t+t^2$. 
Let $c(D)$ be the crossing number of $D$. 
We have 
$$\lim _{t\to 1} X_D(t)=c(D)$$
by definition. 
Hence $X_D(t)$ is a quantization of the crossing number of $D$. 
Let $e$ be an edge of $D$. 
We denote $d(D_b)$ by $d(e)$, where $b$ is a base point on $e$. 
Let $P$ be a projection of a knot with the crossing number $c(P)=n\ge 1$. 
We obtain $2^n$ diagrams $D$ from $P$ by giving over/under information to each double point 
as shown in Figures \ref{s-sum1}, \ref{s-sum2}. 
We call each such diagram $D$ a \textit{state} for $P$. 
Because of the over/under information, states for $P$ have various warping crossing polynomials. 
Then, we consider the \textit{state sum} $Z_P(t)=\sum _D X_D(t)$ of $P$, where $\sum _D$ is the sum for all the states for $P$. 
For example, we have $Z_P(t)=8(1+t)^3$ for the knot projection $P$ with $c(P)=4$ in Figure \ref{s-sum1}. 
We have the following theorem: 

\phantom{x}
\begin{thm}
(i) Let $P$ be a knot projection with $c(P)=n\ge 1$. Then, 
$$Z_P(t)=2n(1+t)^{n-1}.$$

\noindent (ii) Let $D$ be a knot diagram, and $D'$ the diagram obtained from $D$ by a crossing change at a crossing point $p$ of $D$. 
Then, 
\begin{align*}
X_D(t)-tX_{D'}(t)=(1-t)A,
\end{align*}
where $A$ is the sum of $t^{d(e)}$ for all edges $e$ from the under-crossing of $p$ to the over-crossing of $p$. 
\label{mainthm}
\end{thm}
\phantom{x}

\begin{figure}[htbp]
 \begin{center}
  \includegraphics[width=30mm]{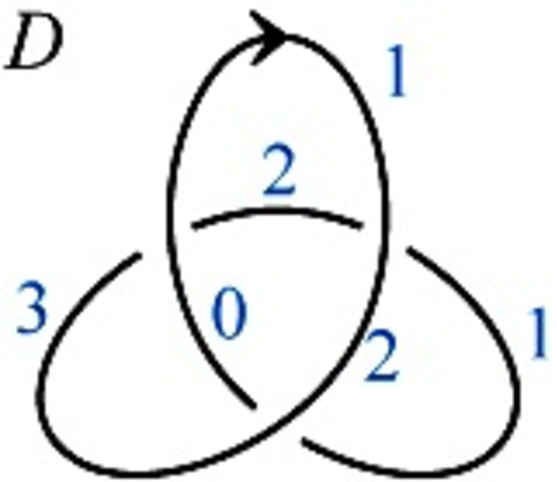}
 \end{center}
 \caption{}
 \label{w-ex}
\end{figure}

\noindent The proof is given in Section \ref{s-state}. 
The \textit{warping polynomial} $W_D(t)$ of a knot diagram $D$ is the sum of $t^{d(e)}$  for all edges $e$ \cite{w-poly}. 
For example, the diagram $D$ in Figure \ref{w-ex} has $W_D(t)=1+2t+2t^2+t^3$. 
We have the following theorem:

\phantom{x}
\begin{thm}
Let $D$ be a knot diagram with $c(D)\ge 1$. 
We have 
$$X_D(t)=\frac{W_D(t)}{1+t}.$$
\label{xw-thm}
\end{thm}
\phantom{x}

\noindent The proof is given in Section \ref{s-wp}. 
The rest of this paper is organized as follows: 
In Section 2, we study a state sum for a plane curve by considering knot diagrams obtained from the plane curve as states. 
In Section 3, we consider properties of the warping crossing polynomial by comparing with the warping polynomial. 
In Section 4, we show that every based plane curve in $\mathbb{R}^2$ has a canonical orientation.

\section{State sum}
\label{s-state}

In this section, we study knot projections by considering the distribution of the states.


\begin{figure}[htbp]
 \begin{center}
  \includegraphics[width=130mm]{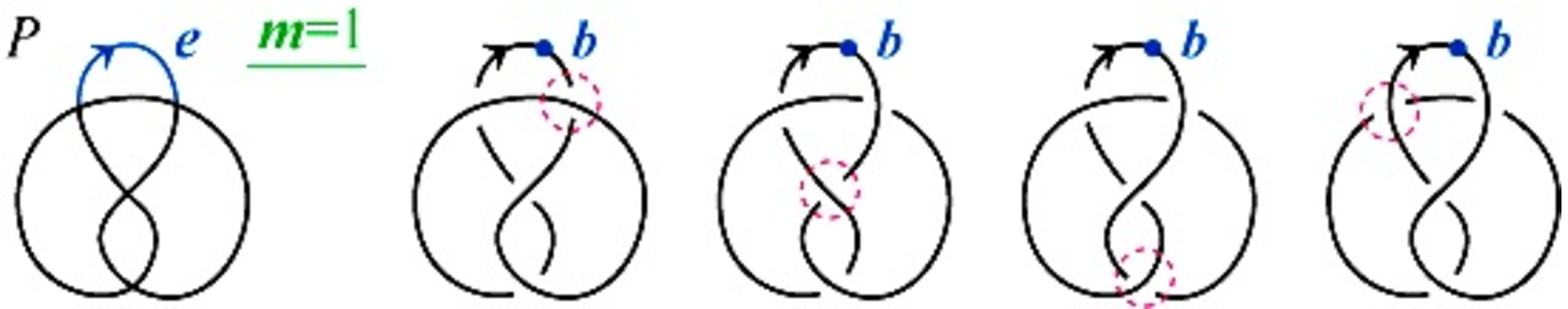}
 \end{center}
 \caption{}
 \label{ss-pf}
\end{figure}

\phantom{x}
\noindent \textit{Proof of Theorem \ref{mainthm}.} 
(i) We show that the sum $\sum _D W_D(t)$ of the warping polynomials $W_D(t)$ for all the states $D$ for $P$ is $2n(1+t)^n$. 
Let $e$ be an edge of $P$, and let $m=1,2,\dots $ or $n$. 
We can give all the double points of $P$ over/under information 
so that $d(e)=m$ in ${}_nC_m$ ways as shown in Figure \ref{ss-pf}. 
Hence 
\begin{align*}
\sum _D W_D(t)&=2n\times {}_nC_0+2n\times {}_nC_1t+2n\times {}_nC_2t^2+\dots +2n\times {}_nC_nt^n\\
                     &=2n(1+t)^n
\end{align*}
because $P$ has $2n$ edges.

\noindent (ii) Let $A$ (resp. $B$) be the sum of $t^{d(e)}$ (resp. $t^{d(e)-1}$) for all edges $e$ 
from the under-crossing (resp. over-crossing) of $p$ to the over-crossing (resp. under-crossing) of $p$. 
By the proof of Lemma 4.4 in \cite{w-poly} and Theorem \ref{xw-thm}, we have $(t+1)X_D(t)=A+tB$ and $(t+1)X_{D'}(t)=tA+B$,  
and therefore we have $X_D(t)-tX_{D'}(t)=(1-t)A$, $X_{D'}(t)-tX_D(t)=(1-t)B$, and $X_D(t)+X_{D'}(t)=A+B$. 
Hence, only the equation $X_D(t)-tX_{D'}(t)=(1-t)A$ is sufficient. 
\hfill$\square$

\phantom{x}

\noindent Let $\mathrm{span}f(t)$ be the span of a polynomial $f(t)$. 
We have the following corollary: 

\phantom{x}
\begin{cor}
Let $D$ and $D'$ be diagrams as above. 
We have 
$$\vert \mathrm{span} X_{D'}(t)-\mathrm{span} X_D(t)\vert \le 2.$$
\end{cor}

\section{Warping crossing polynomial}
\label{s-wp}

In this section, we prove Theorem \ref{xw-thm} and show properties of the warping crossing polynomial. 
We prove Theorem \ref{xw-thm}.

\begin{figure}[htbp]
 \begin{center}
  \includegraphics[width=70mm]{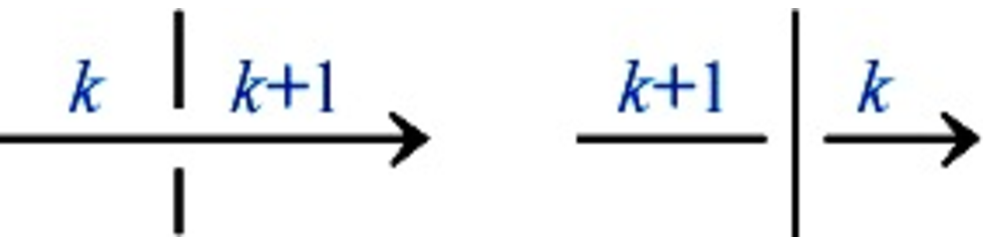}
 \end{center}
 \caption{}
 \label{xw-pf}
\end{figure}

\phantom{x}
\noindent \textit{Proof of Theorem \ref{xw-thm}.} 
If $D$ has $n$ over-crossings shown in the left hand in Figure \ref{xw-pf}, 
then $D$ has also $n$ under-crossings shown in the right hand of Figure \ref{xw-pf}. 
In other words, if there are $n$ edges $e$ such that $d(e)=k$ and the endpoints are over-crossings, 
then there are also $n$ edges $e$ such that $d(e)=k+1$ and the endpoints are under-crossings. 
Since the crossing weight of the crossing point of the left hand of Figure \ref{xw-pf} is $t^k$, 
the sum of $t^{d(e)}$ for all the edges $e$ of $D$ whose endpoints are over-crossings is equal to $\sum _c t^{d(c)}=X_D(t)$, 
and therefore that for all the edges $e$ of $D$ whose endpoints are under-crossings is $\sum _c t^{d(c)+1}=tX_D(t)$. 
Hence $W_D(t)$, which is the sum of $t^{d(e)}$ for all the edges, is $(1+t)X_D(t)$. 
\hfill$\square$

\phantom{x}

\noindent Then, $X_D(t)$ has some properties as $W_D(t)$ has in \cite{w-poly}:

\phantom{x}
\begin{cor}
Let $-D$ be a knot diagram $D$ with the orientation reversed, and $D^*$ the mirror image of $D$. 
We have $X_{-D}(t)=X_{D^*}(t)=t^{n-1}X_D(t^{-1})$, where $n=c(D)$. 
\end{cor}
\phantom{x}

\begin{cor}
A polynomial $f(t)$ is a warping crossing polynomial of a knot diagram $D$ with $c(D)=n\ge 1$ if and only if 
$f(t)=m_0t^d+m_1t^{d+1}+\dots +m_st^{d+s}$, where $m_i=1,2,\dots (i=0,1,\dots ,s)$, $d, s=0,1,\dots $ and 
$m_0+m_1+\dots +m_s=n$. 
\end{cor}
\phantom{x}

\noindent A knot diagram $D$ is an \textit{alternating diagram} if we come to crossing points as an over-crossing and as an under-crossing alternately when we go along $D$. 
A \textit{bridge} in a knot diagram $D$ is a path on $D$ between under-crossings which has no under-crossings and at least one over-crossing in the interior. 
A knot diagram $D$ is a \textit{one-bridge diagram} if $D$ has exactly one bridge. 
The warping crossing polynomial characterizes an alternating diagram and a one-bridge diagram 
as the warping polynomial characterizes in \cite{w-poly}: 

\phantom{x}
\begin{cor}
A knot diagram $D$ with $c(D)=n\ge 1$ is an alternating diagram if and only if $X_D(t)=nt^d$ $(d=0,1,\dots )$. 
\end{cor}
\phantom{x}
\begin{rem}
An alternating diagram $D$ with $c(D)\ge 1$ has constant crossing weights at all the crossing points (see Figure \ref{alt34}). 
\begin{figure}[htbp]
 \begin{center}
  \includegraphics[width=60mm]{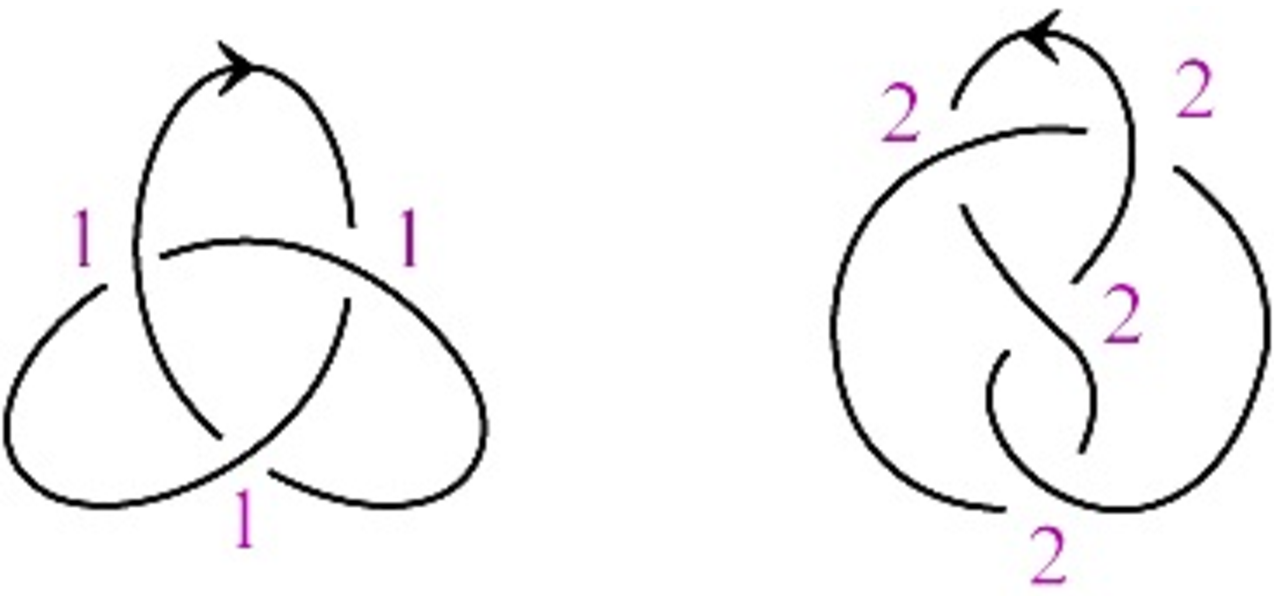}
 \end{center}
 \caption{}
 \label{alt34}
\end{figure}
\end{rem}
\phantom{x}

\begin{cor}
A knot diagram $D$ with $c(D)=n\ge 1$ is a one-bridge diagram if and only if $X_D(t)=1+t+t^2+\dots +t^{n-1}$. 
\end{cor}
\phantom{x}

\begin{rem}
A one-bridge diagram has different crossing weights at all the crossing points (see Figure \ref{one34}). 
\begin{figure}[htbp]
 \begin{center}
  \includegraphics[width=60mm]{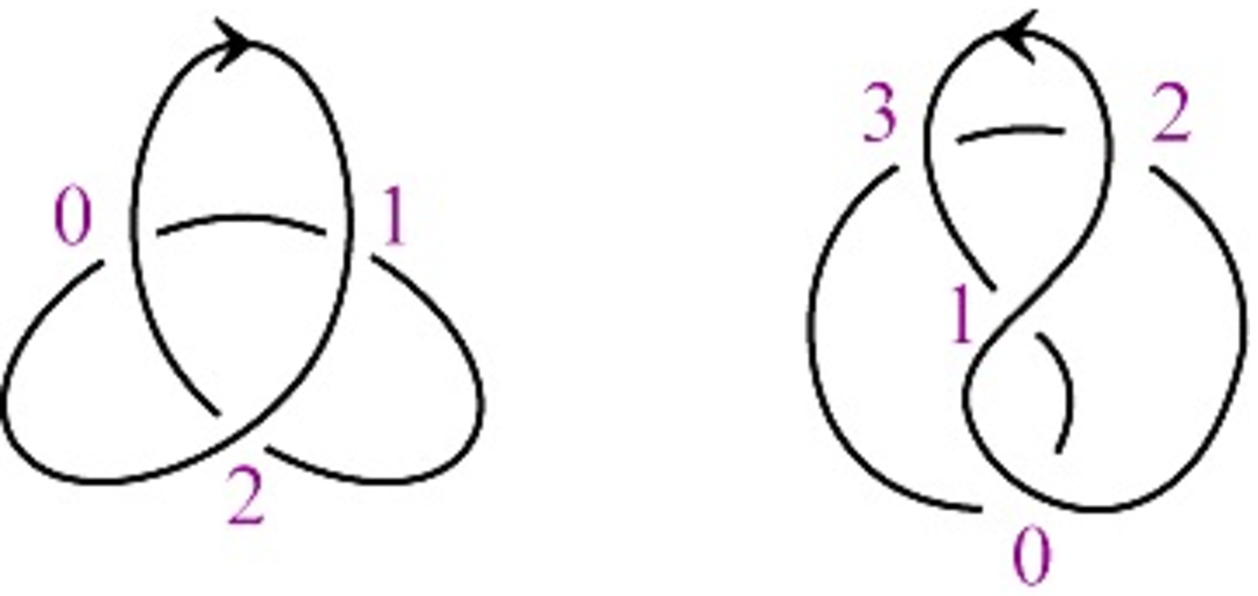}
 \end{center}
 \caption{}
 \label{one34}
\end{figure}
\end{rem}
\phantom{x}

\begin{figure}[htbp]
 \begin{center}
  \includegraphics[width=75mm]{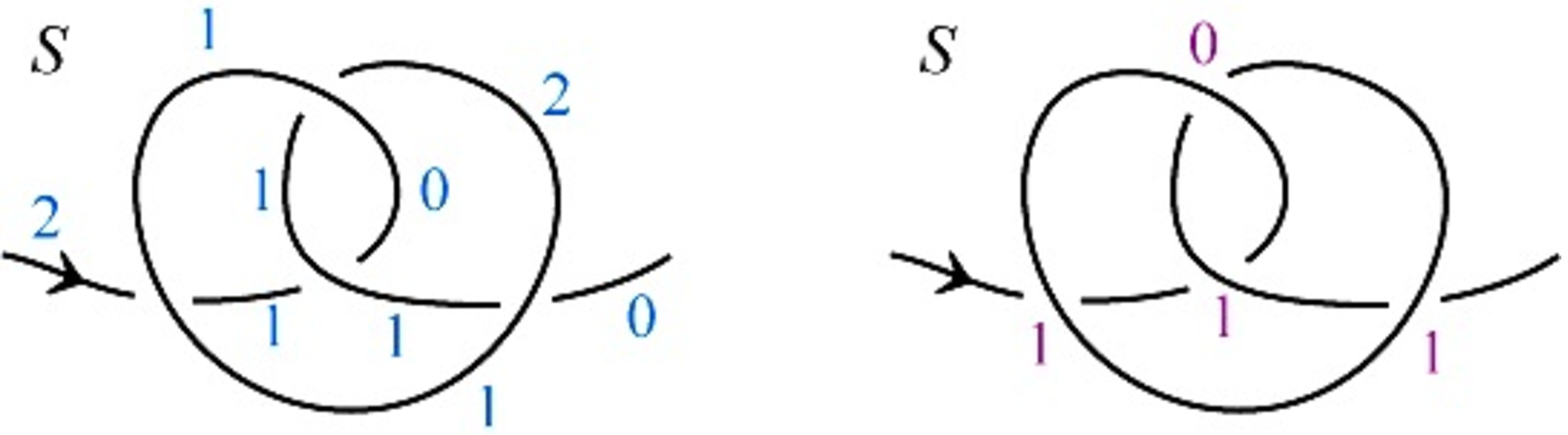}
 \end{center}
 \caption{}
 \label{spatial-a}
\end{figure}

\noindent A \textit{spatial arc diagram} is a diagram of a spatial arc. 
we remark that we can define the warping polynomial $W_S(t)$ and the warping crossing polynomial $X_S(t)$ of a spatial arc diagram $S$. 
For example, we have $W_S(t)=2+5t+2t^2$ and $X_S(t)=1+3t$ for the spatial arc diagram $S$ in Figure \ref{spatial-a}. 

\section{Orientations of plane curves}
\label{s-ori}

In this section we show that we can give each based plane curve on $\mathbb{R}^2$ a canonical orientation by using the warping degrees. 
We first review the warping degree of a (non-based) knot diagram. 
The \textit{warping degree} $d(D)$ of an oriented knot diagram $D$ is the minimal warping degree $d(D_b)$ of $D_b$ for all base points $b$ of $D$ \cite{kawauchi-l}. 
The following theorem is shown in \cite{shimizu1}: 

\phantom{x}
\begin{thm}
Let $D$ be an oriented knot diagram with $c(D)\ge 1$. 
We have 
$$d(D)+d(-D)+1\le c(D).$$
Further, the equality holds if and only if $D$ is an alternating diagram. 
\end{thm}
\phantom{x}

\noindent We have the following corollary: 

\phantom{x}
\begin{cor}
Let $D$ be an oriented alternating knot diagram with non-zero even crossings. 
Then, $X_D(t)\ne X_{-D}(t)$. 

\phantom{x}
\begin{proof}
We have $d(D)+d(-D)=c(D)-1$ because $D$ is alternating. 
Since $c(D)$ is even, the value $d(D)+d(-D)$ is odd. 
Hence the crossing weights of the crossing points of $-D$ are different from that of $D$. 
\end{proof}
\label{d-d}
\end{cor}
\phantom{x}

\noindent Now we discuss the orientations of plane curves. 
We have the following theorem: 

\phantom{x}
\begin{thm}
(1) Every closed transversely intersected curve $C$ with even crossing points on $\mathbb{R}^2$ has two independent canonical orientations. 

\noindent (2) Every based closed transversely intersected curve $C_b$ with odd crossing points on $\mathbb{R}^2$ has two independent canonical orientations. 

\phantom{x}
\begin{proof} 
(1) We give $C$ with even crossing points one orientation in the following order: 
First, we explain how to obtain an alternating diagram uniquely from $C$. 
After that, we give the alternating diagram the canonical orientation. 
Apply $C$ the checkerboard coloring such that the outer region is colored white. 
Then we obtain an alternating diagram $D$ uniquely by giving each double point over/under information as shown in Figure \ref{ori-8}. 
If  $C$ has no crossing point, then we consider as $D$ a knot diagram obtained from $C$ by taking connected sums with two positive one-crossing knot diagrams in the black region. 
\begin{figure}[htbp]
 \begin{center}
  \includegraphics[width=70mm]{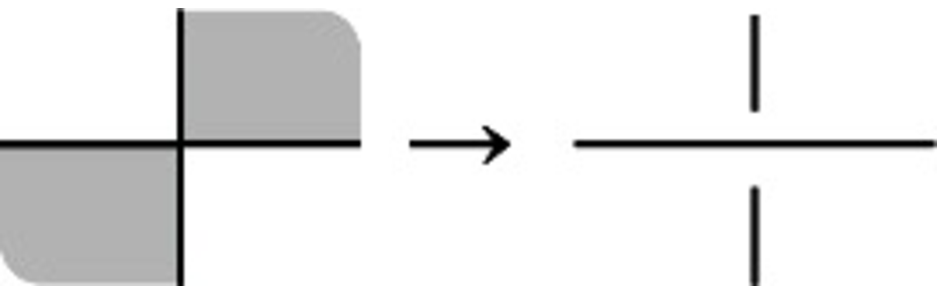}
 \end{center}
 \caption{}
 \label{ori-8}
\end{figure}
Since $c(D)$ is non-zero even, we can give $D$ the orientation uniquely so that $d(D)<d(-D)$ by the proof of Corollary \ref{d-d}. 
By projection, $C$ is also oriented. \\
The other orientation is given by using the rotation number. 
The rotation number $\mathrm{rot}(E)$ of an oriented closed curve $E$ on $\mathbb{R}^2$ is $l_+(E)-l_-(E)$, where $l_+(E)$ (resp. $l_-(E)$) is the number of circles with counter-clockwise (resp. clockwise) orientations obtained by splicing $E$ at all the crossing points. 
Note that $\mathrm{rot}(E)$ is odd if $c(E)$ is even, and that we have $\mathrm{rot}(-E)=-\mathrm{rot}(E)$. 
Hence for non-oriented closed curve $C$, we can give $C$ the orientation uniquely so that $\mathrm{rot}(C)$ is positive. 
Figure \ref{ori-even} shows that these two orientations of $C$ are independent. 
\begin{figure}[htbp]
 \begin{center}
  \includegraphics[width=110mm]{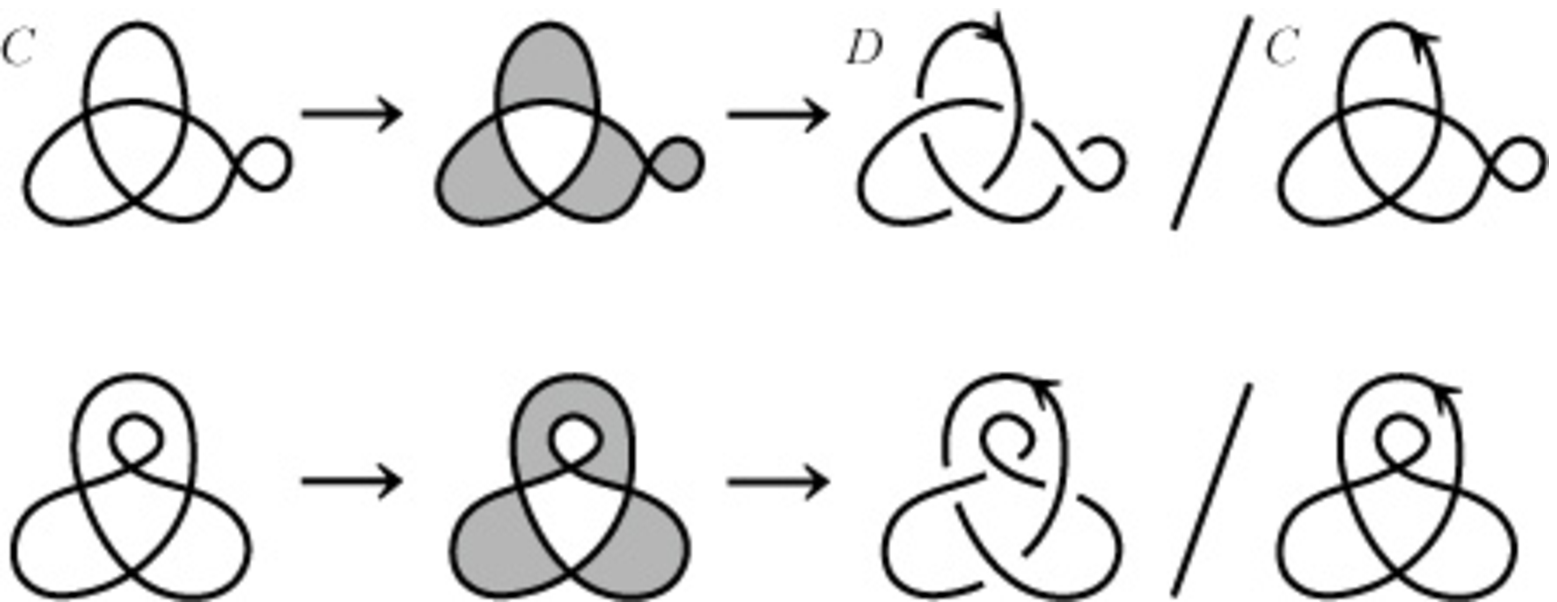}
 \end{center}
 \caption{}
 \label{ori-even}
\end{figure}

\noindent (2) We give $C_b$ with odd crossing points one orientation as follows: 
We apply $C_b$ the checkerboard coloring as above and we obtain the alternating diagram $D_b$. 
Apply the connected sum of a knot diagram with exactly one positive crossing to the edge with the base point $b$ in the black region, and we obtain $D'$ (see Figure \ref{ori-odd}). 
Since $D'$ is alternating and with even crossings, $D'$ has the canonical orientation. 
Therefore $D_b$ and $C_b$ are also oriented. \\
The other orientation of $C_b$ with the checkerboard coloring is the orientation such that at the base point $b$ the black region lies on the right. 
Figure \ref{ori-odd} shows that these two orientations of $C_b$ are independent. 
\begin{figure}[htbp]
 \begin{center}
  \includegraphics[width=110mm]{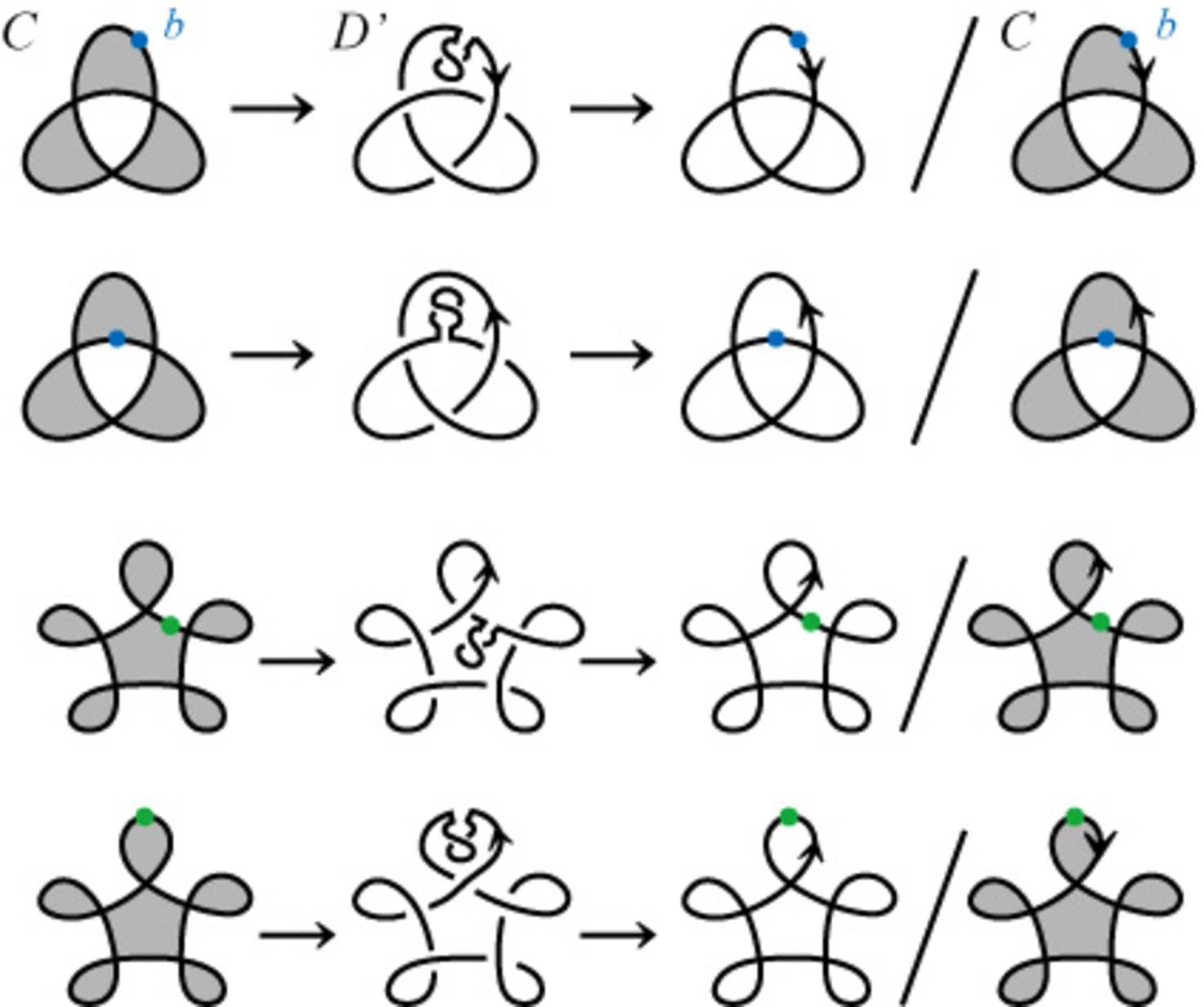}
 \end{center}
 \caption{}
 \label{ori-odd}
\end{figure}
\end{proof}
\label{ori-thm}
\end{thm}
\phantom{x}

\noindent 
The first orientations in (1) and (2) were given in the first version of this paper and the second orientations in (1) and (2) were suggested later by K. Taniyama and V. Turaev, respectively. 
It is an interesting question to explain a difference between the two independent orientations in (1) and (2). 
We have the following corollary: 

\phantom{x}
\begin{cor}
(1) For every oriented closed transversely intersected curve $C$ with even crossing points on $\mathbb{R}^2$, there is no orientation-preserving homeomorphism from $\mathbb{R}^2$ to $\mathbb{R}^2$ sending $C$ to $-C$.

\noindent (2) For every based oriented closed transversely intersected curve $C_b$ with odd crossing points on $\mathbb{R}^2$, there is no orientation-preserving, base-point-preserving homeomorphism from $\mathbb{R}^2$ to $\mathbb{R}^2$ sending $C_b$ to $-C_b$. 
\label{ori-cor}
\end{cor}

\phantom{x}
\begin{rem}
If $C$ with odd crossing points is non-based, the corollary above does not hold (see Figure \ref{cor-hom}). 
\begin{figure}[htbp]
 \begin{center}
  \includegraphics[width=80mm]{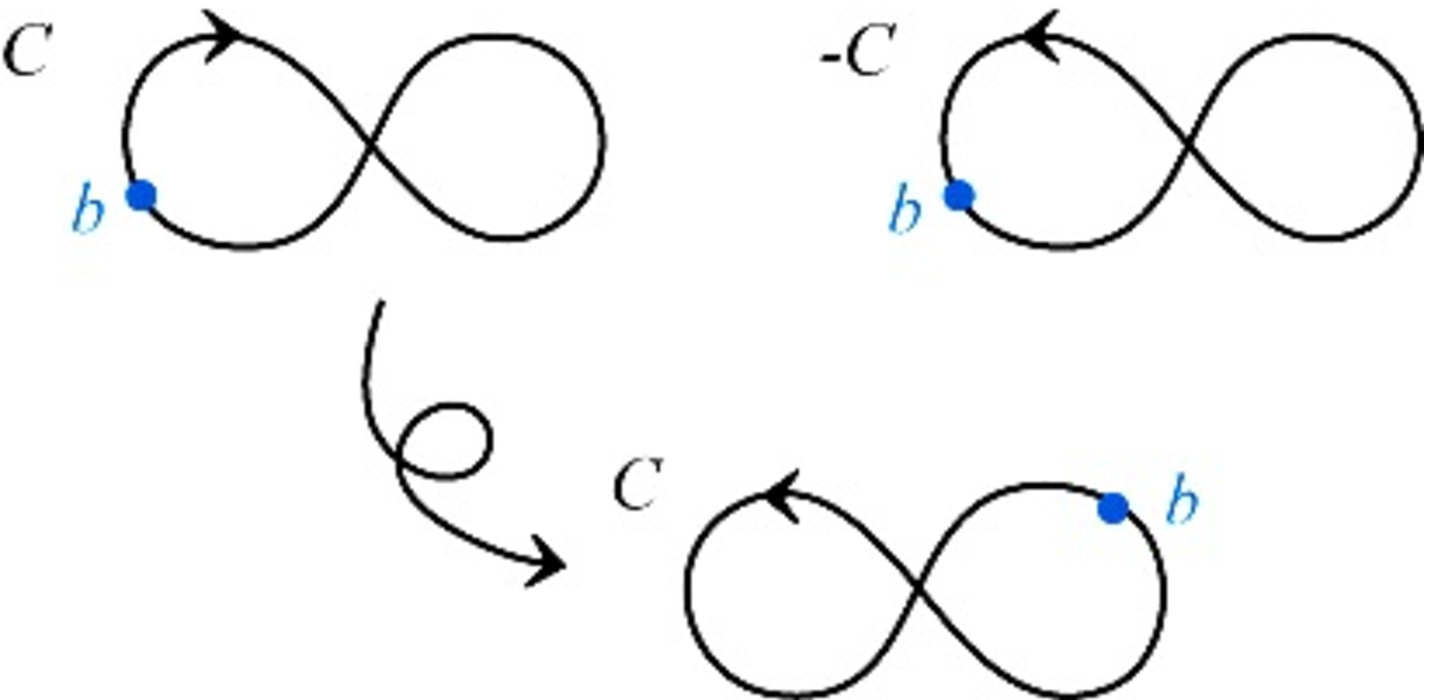}
 \end{center}
 \caption{}
 \label{cor-hom}
\end{figure}
\end{rem}
\phantom{x}

\begin{rem}
Theorem \ref{ori-thm} (2) and Corollary \ref{ori-cor} (2) hold on $S^2$ in place of $R^2$. 
In fact, the number of regions divided by $C$ is $c(C)+2$ which is odd. 
Hence we can take the unique checkerboard coloring so that the number of black regions is greater than the number of white regions.
\end{rem}
\phantom{x}

\section*{Acknowledgments}
The authors would like to thank Professor K. Taniyama and Professor V. Turaev  for helpful suggestions on Theorem \ref{ori-thm}. 
The first author is supported by Grant-in-Aid for Scientific Research (A) (21244005). 
The second author is partly supported by JSPS Research Fellowships for Young Scientists.

\maketitle


\begin{thebibliography}{9}
\bibitem{fujimura}S.~Fujimura: \textit{On the ascending number of knots}, thesis, Hiroshima University, 1988. 
\bibitem{fung}T.~S.~Fung: \textit{Immersions in knot theory}, a dissertation, Columbia University, 1996. 
\bibitem{kawauchi-s}A.~Kawauchi: \textit{A survey of knot theory}, Birkhauser, (1996). 
\bibitem{kawauchi-l}A.~Kawauchi: \textit{Lectures on knot theory} (in Japanese), Kyoritsu shuppan Co. Ltd, 2007.
\bibitem{kawauchi-c}A.~Kawauchi: \textit{On a complexity of a spatial graph}, in: Knots and softmatter physics, Topology of polymers and related topics in physics, mathematics and biology, Bussei Kenkyu 92-1 (2009-4), 16--19.
\bibitem{LM}W.~B.~R.~Lickorish and K.~C.~Millett: \textit{A polynomial invariant of oriented links}, Topology {\bf 26} (1987), 107--141. 
\bibitem{okuda}M.~Okuda: \textit{A determination of the ascending number of some knots}, thesis, Hiroshima University, 1988. 
\bibitem{ozawa}M. Ozawa: \textit{Ascending number of knots and links}, J. Knot Theory Ramifications {\bf 19} (2010), 15--25. 
\bibitem{polyak}M.~Polyak: \textit{Minimal generating sets of Reidemeister moves}, arXiv:0908.3127v3. 
\bibitem{shimizu1}A.~Shimizu: \textit{The warping degree of a knot diagram}, J. Knot Theory Ramifications, {\bf 19} (2010), 849--857. 
\bibitem{shimizu2}A.~Shimizu: \textit{The warping degree of a link diagram}, Osaka J. Math., {\bf 48} (2011), 209--231. 
\bibitem{w-poly}A.~Shimizu: \textit{The warping polynomial of  a knot diagram}, arXiv:1109.5898v1.
\end{thebibliography}
\end{document}